\documentclass[12pt,reqno]{amsart}

\usepackage[english]{babel}
\usepackage[utf8]{inputenc}
\usepackage{amsmath}
\usepackage{amssymb}
\usepackage{amsthm}
\usepackage{enumerate}
\usepackage{graphicx}
\usepackage{fullpage}
\usepackage{url}

\newtheorem{thm}{Theorem}[section]
\newtheorem{lemma}[thm]{Lemma}

\newtheorem{conjecture}[thm]{Conjecture}

\newcommand{\R}{\mathbb{R}}

\newcommand{\mP}{\mathcal{P}}

\newcommand{\mR}{\mathcal{R}}

\newcommand{\sm}{\setminus}

\title{Maximum number of colourings. II. 5-chromatic graphs}

\author{Fiachra Knox}
\address{Department of Mathematics\\
      Simon Fraser University\\
      Burnaby, BC \ V5A 1S6}
\email{fiachraknox@hotmail.com}
\thanks{F.K.~was supported by a PIMS Postdoctoral Fellowship.}

\author{Bojan Mohar}
\address{Department of Mathematics\\
      Simon Fraser University\\
      Burnaby, BC \ V5A 1S6}
\email{mohar@sfu.ca}
\thanks{B.M.~was supported in part by the NSERC Discovery Grant R611450 (Canada),
by the Canada Research Chairs program, and by the Research Project J1-8130 of ARRS (Slovenia).}
\thanks{On leave from IMFM \& FMF, Department of Mathematics, University of Ljubljana.}%

\date{\today}

\begin{document}

\begin{abstract}
In 1971, Tomescu conjectured [Le nombre des graphes connexes $k$-chromatiques minimaux aux sommets \'{e}tiquet\'{e}s, C. R. Acad. Sci. Paris 273 (1971), 1124--1126] that every connected graph $G$ on $n$ vertices with $\chi(G) = k \geq 4$ has at most $k!(k-1)^{n-k}$ $k$-colourings, where equality holds if and only if the graph is formed from $K_k$ by repeatedly adding leaves.
In this note we prove (a strengthening of) the conjecture of Tomescu when $k=5$.
\end{abstract}

\maketitle

\section{Introduction}

Let $x$ be a positive integer. By an \emph{$x$-colouring} we mean a function $f:V(G)\to\{1,\dots,x\}$ such that $f(u)\ne f(v)$ whenever $uv\in E(G)$. Note that permuting the colours used in a colouring gives a different colouring. The \emph{chromatic polynomial\/} $P_G(x)$ is the polynomial of degree $n=|V(G)|$ whose value $P_G(x)$ is equal to the number of $x$-colourings of $G$ for every positive integer $x$.  Throughout the paper we will use the related indeterminate $y=x-1$ and the \emph{shifted chromatic polynomial}:
$$
   Q_G(y) = P_G(y+1).
$$

Very basic questions about chromatic polynomials remain unresolved and poorly understood.  We refer to part I \cite{KM1} for history and motivation related to using the chromatic polynomial.  In this paper we continue the work on
maximizing the number of colourings among all connected graphs of given order, with the goal to prove a conjecture of Tomescu \cite{To71} dating back to 1971.
The conjecture is motivated by the following easy fact.
For every $k\ge1$ and every integer $x\ge k$, every connected $n$-vertex graph containing a clique of order $k$ has at most
\begin{equation}
   x\cdot(x-2)(x-3)\cdots (x-k+1)\cdot (x-1)^{n-k+1} = x^{\underline{k}}\,(x-1)^{n-k}
\label{eq:1}
\end{equation}
$x$-colourings, where $x^{\underline{k}} = x(x-1)(x-2)\cdots (x-k+1)$ is the falling factorial. This bound is attained for every $x$ if $G$ can be obtained from the $k$-clique $K_k$ by growing an arbitrary tree from each vertex of the clique. In 1971, Tomescu \cite{To71} conjectured that (\ref{eq:1}) is an upper bound for the number of $k$-colourings of any connected $k$-chromatic graph, whether it contains a $k$-clique or not, as long as $k\ge4$:

\begin{conjecture}[Tomescu, 1971]
Let\/ $G$ be a connected $k$-chromatic graph with $k\ge4$. Then $G$ has at most
\begin{equation}
   k!(k-1)^{|G|-k}
\label{eq:2}
\end{equation}
$k$-colourings. Moreover, the extremal graphs are precisely the graphs obtained from $K_k$ by adding trees rooted at each vertex of the clique.
\end{conjecture}

The requirement that $k\ne3$ is necessary since odd cycles of length at least $5$ (and graphs formed from them by adding trees rooted at their vertices) have more colourings than specified by (\ref{eq:1}); see \cite{To_DM72} for more details; and in the case of bipartite graphs ($k=2$), any connected bipartite graph attains the bound.

Tomescu proved \cite{To_JGT90} that all 4-chromatic planar graphs satisfy his conjecture. For the same class of graphs, he proved a stronger conclusion for the number or $x$-colourings (for every $x\ge4$), where the bound of the conjecture is replaced by (\ref{eq:1}). Apart from this achievement, only sporadic results are known \cite{BE15,BEL16}. We refer to \cite[Chapter~15]{DKT_Book_05} for additional overview of the results in this area.

In \cite{KM1} we proved an extended version of Tomescu conjecture for $k=4$.

\begin{thm}[\cite{KM1}]
\label{thm:4chrom}
Let $G$ be a connected $4$-chromatic graph and $y\ge 3$ be an integer. Then
\begin{equation}
   Q_G(y) \le (y+1)y^{n-3}(y-1)(y-2).
\label{eq:3}
\end{equation}
Moreover, equality holds for some integer $y\ge 3$ if and only if $G$ can be obtained from $K_4$ by adding a tree on each vertex of $K_4$ (in which case equality holds for every $y \in \R$).
\end{thm}

In this note we settle the case when $k=5$, again in its extended version for colourings with any number of colours.

\begin{thm}
\label{thm:main}
Let $G$ be a connected $5$-chromatic graph and $y\ge 4$ be an integer. Then
\begin{equation}
   Q_G(y) \le (y+1)y^{n-4}(y-1)(y-2)(y-3).
\label{eq:3a}
\end{equation}
Moreover, equality holds for some integer $y\ge 4$ if and only if $G$ can be obtained from $K_5$ by adding a tree on each vertex of $K_5$ (in which case equality holds for every $y \in \R$).
\end{thm}

\section{Preliminaries}

We will use standard graph theory terminology and notation as used by Diestel \cite{Diestel2010} or Bondy and Murty \cite{BondyMurty2008}. In particular, we use $n=|G|=|V(G)|$ to denote the order of $G$. The minimum vertex degree of $G$ is denoted by $\delta(G)$. By $N(v)$ we denote the set of neighbours of a vertex $v$, and we use $\chi(G)$ for the chromatic number. We say $G$ is \emph{$k$-chromatic} if $\chi(G)=k$. The graph is \emph{vertex-critical} (\emph{edge-critical}) if the removal of any vertex (edge) decreases the chromatic number. We will frequently use the fact that identifying non-adjacent vertices of a graph $G$ results in a graph $G'$ with $\chi(G') \geq \chi(G)$.
For a vertex-set $U\subseteq V(G)$, $G[U]$ is the subgraph of $G$ induced by $U$.

If $t$ is a positive integer, we let $[t] = \{1,2,\dots,t\}$.

\begin{lemma} \label{lem:single vx bounds}
Let $G$ be a graph, let $v \in V(G)$ have degree $d$ and let $G' = G - v$.
Suppose that for every $\ell \ge 1$ and every graph $H$ on $|G'|-\ell$ vertices that is formed from $G'$ by repeatedly identifying pairs of nonadjacent vertices satisfies
$$Q_H(y) \le y^{-\ell} B(y).$$
For $r \in [d-1]$, let $N_r=N_r(v)$ be the number of partitions of $N(v)$ into $r$ non-empty independent sets. Then
$$Q_G(y) \leq (y - r) Q_{G'}(y) + N_{r} y^{-d+r} B(y).$$
\end{lemma}

\proof For each $i \in [d]$, let $Q_i$ be the number of $(y+1)$-colourings of $G'$ which take exactly $i$ colours on $N(v)$.
Then
\begin{equation} \label{eq:svb1}
Q_G(y) = \sum_{i = 1}^{d} (y - i + 1)Q_i
\leq (y-r) Q_{G'}(y) + \sum_{j = 0}^{r-1} (j+1) Q_{r-j},
\end{equation}
since $\sum_{i = 1}^d Q_i = Q_{G'}(y)$.
(Note that if $y < r$, then $Q_i = 0$ for every $i > r$ and so equality holds in (\ref{eq:svb1}).)

Let $\Omega_r$ be the set of partitions $\mP$ of $N(v)$ into exactly $r$ non-empty independent sets, and let $\Omega = \bigcup_{r = 1}^d \Omega_r$.
For each $\mP \in \Omega_r$, let $G_{\mP}$ be the graph formed from $G'$ by identifying the vertices in each part of $\mP$
to a single vertex,
and let $G_{\mP}^*$ be the graph formed from $G_{\mP}$ by adding edges between every pair of non-adjacent identified vertices.
Given a partition $\mR$ of $N(v)$ into non-empty independent sets, we write $\mP \geq \mR$ if $\mP$ refines $\mR$.
For brevity we write $Q_{\mP}(y)$ and $Q_{\mP}^*(y)$ for $Q_{G_{\mP}}(y)$ and $Q_{G_{\mP}^*}(y)$, respectively. Now
\begin{equation} \label{eq:svb2}
\sum_{j=0}^{r-1} (j+1) Q_{r - j} \leq
\sum_{\mP \in \Omega_r} \sum_{\mP \geq \mR \in \Omega} Q_{\mR}^*(y) =
\sum_{\mP \in \Omega_r} Q_{\mP}(y) \leq
N_r y^{-d + r} B(y),
\end{equation}
where the first inequality holds since any partition $\mR \in \Omega_{r - j}$ can be refined into a partition into $r$ parts in at least $j+1$ ways, and hence the term $Q_{\mR}^*(y)$ appears at least $j+1$ times in the sum $\sum_{\mP \in \Omega_r} \sum_{\mP \geq \mR \in \Omega} Q_{\mR}^*(y)$.
The desired inequality now follows from (\ref{eq:svb1}) and (\ref{eq:svb2}).
\endproof

\section{Proof of Theorem \ref{thm:main}}

We define a partial ordering $\ll_k$ on polynomials in $y$ by $P_2 \ll_k P_1$ when every coefficient of $W(z) = (P_1 - P_2)(z + k)$ is non-negative.
Note that this implies that $P_1(y)\ge P_2(y)$ for every $y\ge k$.
We write $\ll$ for $\ll_4$, as in the majority of cases we will take $k = 4$.

\begin{lemma} \label{lem:comp1}
Let $G$ be a $5$-critical graph.
If there exists $S \subseteq V(G)$ such that $|S| \geq 5$ and $\chi(G - S) \ge 4$, then
$$Q_G(y) \leq (y+1)y^{n-4}(y-1)(y-2)(y-3).$$
\end{lemma}

\proof We may assume that $G - S$ is connected; indeed, if $G - S$ is not connected then we can add to $S$ the vertices of each component but one (chosen to be $4$-chromatic).
We may also assume that $|S| = 5$, noting that we can reduce the size of $S$ (if necessary) by repeatedly removing vertices with a neighbour outside $S$,.

For each vertex $v \in S$, let $d'(v) = 4 - |N(v) \cap S|$ and note that $v$ has at least $d'(v)$ neighbours outside $S$.
Let $S'$ be the subset of $S$ consisting of vertices $v$ with $d'(v) \geq 2$ (i.e., those vertices with $|N(v) \cap S| \leq 2$).
For any $T \subseteq S'$, Theorem~\ref{thm:4chrom} implies that the number of colourings $F_T$ of $G - S$ in which $N(v) \sm S$ is monochromatic for every $v \in T$ is at most
$(y+1)y^{n-\Delta'(T)-7}(y-1)(y-2)$, where $\Delta'(T) = \max \{d'(v) \mid {v \in T}\} \cup \{1\}$, since every such colouring can be derived from a colouring of a graph formed from $G - S$ by identifying an independent set of $\Delta'(T)$ vertices (note that identifying an independent set preserves connectivity and does not decrease the chromatic number).

Further, for any $T \subseteq S'$ we can compute a polynomial  upper bound $E_T$ on the number of extensions of a colouring of $G - S$ to $S$, given that $N(v) \sm S$ is not monochromatic for any $v \in S' \sm T$.
We do this by assigning to each vertex of $S$ a set of $0$, $1$ or $2$ forbidden colours (as appropriate), computing an upper bound on the number of extensions by deletion-contraction (where an isolated vertex with $r$ forbidden colours is given the upper bound $y - r + 1$, and contracted edges produce a vertex whose set of forbidden colours is the union of the sets of forbidden colours of the endvertices), and selecting those of the resulting polynomials (in $y$) which are maximal under $\ll$.
Even though $\ll$ is not a total ordering, in this case it turns out that each of these maximal polynomials is unique.

For each $T \subseteq S'$, let $F'_T$ be the number of colourings of $G - S$ in which $N(v) \sm S$ is monochromatic for every $v \in T$, but not for any $v \in S' \sm T$.
Now we have $Q_G(y) \leq \sum_{T \subseteq S'} E_T F'_T$.
Observe that $F_T = \sum_{T \subseteq T^+ \subseteq S'} F'_{T^+}$.
Let $E'_T = \sum_{T^- \subseteq T} (-1)^{|T \sm T^-|} E_T$ for each $T \subseteq S'$, so that $E_T = \sum_{T^- \subseteq T} E'_{T^-}$ for each $T \subseteq S'$.
Then
$$\sum_{T \subseteq S'} E_T F'_T = \sum_{T^- \subseteq T^+ \subseteq S'} E'_{T^-} F'_{T^+} = \sum_{T \subseteq S'} E'_T F_T.$$
In general the polynomials $E'_T$ may be positive or negative.
Fortunately, each $E'_T$ is either positive for every integer $y \geq 4$, or negative for every integer $y \geq 4$.
Our computer program has verified this using the relation $\ll_k$ (for appropriate values of $k$) along with individual checks for small values of $y$.
Let $\Omega_+$(resp. $\Omega_-$) be the family of sets $T \subseteq S'$ such that $E'_T$ is positive (resp. negative) for every integer $y \geq 4$.
Then
\begin{eqnarray*}
  Q_G(y) &\leq& \sum_{T \in \Omega_+} E'_T (y+1)y^{n-\Delta'(T)-7}(y-1)(y-2) \\
    &=& (y+1)y^{n-11}(y-1)(y-2) \sum_{T \in \Omega_+} E'_T y^{4-\Delta'(T)}.
\end{eqnarray*}
Observe that $S'$, $\Omega_+$ and $\Omega_-$ depend only on $G[S]$.
Further, $E'_T$ and $\Delta'(T)$ depend only on $G[S]$ and $T$ for each $T \subseteq S'$.
For each of the $34$ (unlabelled) graphs on $5$ vertices, we have computed the polynomial $R=\sum_{T \subseteq S'} E'_T y^{4-\Delta'(T)}$.
The results are given in Table~\ref{tab:1}\footnote{The computer code used can be obtained from the authors.}.
In each case, the resulting polynomial satisfies $R \ll y^7 (y-3)$ (see Table~\ref{tab:2}), and this proves the lemma. \endproof

%The proofs in this paper use Lemma \ref{lem:comp1} only for the case when $S$ contains an independent set of size 4.
%There are just 5 such cases to consider. The values of the corresponding polynomial $Q$ are given in Table \ref{tab:1}.

\begin{lemma} \label{lem:comp2} Let $G$ be a $5$-critical graph.
If $G$ has no clique or independent set of size at least $4$, then $Q_G(y) \ll (y+1)y^{n-4}(y-1)(y-2)(y-3)$.
\end{lemma}

\proof Since $G$ is $5$-critical, removing any vertex leaves a $4$-chromatic graph of order at most $12$, since each colour class has size at most $3$.
Hence $G$ has order at most $13$.
The $5$-edge-critical graphs on at most 12 vertices have already been listed by Royle~\cite{RoyleLink};
for each one\footnote{All together there are 151948 such graphs.}, we have computed its chromatic polynomial $Q$ and verified that $Q \ll (y+1)y^{n-4}(y-1)(y-2)(y-3)$.

It remains only to check graphs on exactly $13$ vertices.
Using Brendan McKay's program \verb nauty_geng, \ one can list all of the Ramsey$(4, 4)$ graphs on $13$ vertices.
In fact these were already listed by McKay~\cite{McKayLink}.
We have tested each such graph for $5$-edge-criticality.
Only $525$ such graphs are $5$-edge-critical.
Again, for each such graph we have computed its chromatic polynomial $Q$.
As it turns out in all cases, we have $Q \ll (y+1)y^{n-4}(y-1)(y-2)(y-3)$. \endproof

\proof[Proof of Theorem \ref{thm:main}]
It is easy to see that we may assume that $G$ is 5-critical.
(The reader may also check \cite{KM1} for details.)
The list of all 5-critical graphs with at most $8$ vertices is available from \cite{RoyleLink}. Their chromatic polynomials satisfy the theorem. This was checked by computer.
Hereafter we may therefore assume that $|G| \geq 9$.
By Lemma~\ref{lem:comp2} we may assume that $G$ has either a clique or an independent set of size $4$.
If $G$ has a clique of size $4$ then we apply Lemma~\ref{lem:comp1} with $S$ being the vertices outside the clique.
So we may assume that $G$ has an independent set $S$ of size $4$.
If $G - S$ is not $4$-critical then we can add a vertex to $S$ so that $G - S$ is still $4$-chromatic, and apply Lemma~\ref{lem:comp1}.
So we may assume that $G- S$ is $4$-critical.

Hence $G - S$ is connected and $4$-chromatic.
By Theorem~\ref{thm:4chrom} we have $Q_{G-S}(y) \leq (y+1)y^{n-7}(y-1)(y-2)$.
Label the vertices of $S$ as $v_1, v_2, v_3, v_4$.
Let $G_i = G - \{v_{i+1}, \ldots, v_4\}$ for $i = 0, 1, 2, 3, 4$.
We may assume that $d(v_i) = 4$ for each $i \in [4]$ (otherwise, we select a subset of $N(v_i)$ and use this in place of $N(v_i)$).
Note that there are $7$ partitions of $N(v_i)$ into two non-empty sets.
For brevity we write $Q_0(y) = (y+1)y^{n-9}(y-1)(y-2)$.
We apply Lemma~\ref{lem:single vx bounds} with $G = G_1$, $v = v_1$, $B(y) = y^2 Q_0(y)$ and $r = 2$ to obtain
$$Q_{G_1}(y) \leq (y-2)y^2 Q_0(y) + 7Q_0(y) = (y^3 - 2y^2 + 7)Q_0(y).$$
We again apply Lemma~\ref{lem:single vx bounds}, this time with with $G = G_2$, $v = v_2$, $B(y) = y^3Q_0(y)$ and $r = 2$ to obtain
$$Q_{G_2}(y) \leq (y-2)(y^3 - 2y^2 + 7)Q_0(y) + 7yQ_0(y) = (y^4 - 4y^3 + 4y^2 + 14y - 14)Q_0(y).$$
We again apply Lemma~\ref{lem:single vx bounds}, this time with with $G = G_3$, $v = v_3$, $B(y) = y^4Q_0(y)$ and $r = 1$ to obtain
\begin{align*}
Q_{G_3}(y) &\leq (y-1)(y^4 - 4y^3 + 4y^2 + 14y - 14)Q_0(y) + yQ_0(y) \\
&= (y^5 - 5y^4 + 8y^3 + 10y^2 - 27y + 14)Q_0(y).
\end{align*}
Finally, we apply Lemma~\ref{lem:single vx bounds}, this time with with $G = G_4$, $v = v_4$, $B(y) = y^5Q_0(y)$ and $r = 1$ to obtain
\begin{align*}
Q_{G_4}(y) &\leq (y-1)(y^5 - 5y^4 + 8y^3 + 10y^2 - 27y + 14)Q_0(y) + y^2Q_0(y) \\
&= (y^6 - 6y^5 + 13y^4 + 2y^3 - 36y^2 + 41y - 14)Q_0(y) \\
&\leq (y^6 - 3y^5)Q_0(y),
\end{align*}
as desired, where the last inequality holds since $y^6 - 6y^5 + 13y^4 + 2y^3 - 36y^2 + 41y - 14 \ll y^6 - 3y^5$.
\endproof

\section{Conclusion}

The proof of the main theorem of this note relies on the 4-chromatic case proved in \cite{KM1}. The main auxiliary Lemmas \ref{lem:comp1} and \ref{lem:comp2} rely on extensive case analysis and use of computer. The main message here is that the same method can be used for larger values of $k$. We see no difficulties of applying it for $k=6$ and possibly for a few additional values.  In a forthcoming paper \cite{KM3} we use the results from \cite{KM1} and from this paper as a basis of induction to tackle the general case.

\section{Acknowledgements}
We thank Jernej Azarija for his significant help and advice in writing the computer programs used to prove Lemmas~\ref{lem:comp1} and~\ref{lem:comp2}.

\bibliographystyle{abbrv}

\bibliography{Tomescu_biblio}

\begin{thebibliography}{10}

\bibitem{BondyMurty2008}
J.~A. Bondy and U.~S.~R. Murty.
\newblock {\em Graph theory}, volume 244 of {\em Graduate Texts in
  Mathematics}.
\newblock Springer, New York, 2008.

\bibitem{BE15}
J.~Brown and A.~Erey.
\newblock New bounds for chromatic polynomials and chromatic roots.
\newblock {\em Discrete Math.}, 338(11):1938--1946, 2015.

\bibitem{BEL16}
J.~I. Brown, A.~Erey, and J.~Li.
\newblock Extremal restraints for graph colourings.
\newblock {\em J. Combin. Math. Combin. Comput.}, 93:297--304, 2015.

\bibitem{Diestel2010}
R.~Diestel.
\newblock {\em Graph theory}, volume 173 of {\em Graduate Texts in
  Mathematics}.
\newblock Springer, Heidelberg, fourth edition, 2010.

\bibitem{DKT_Book_05}
F.~M. Dong, K.~M. Koh, and K.~L. Teo.
\newblock {\em Chromatic polynomials and chromaticity of graphs}.
\newblock World Scientific Publishing Co. Pte. Ltd., Hackensack, NJ, 2005.

\bibitem{KM1}
F.~Knox and B.~Mohar.
\newblock Maximum number of colourings. {I}. 4-chromatic graphs.
\newblock Submitted.

\bibitem{KM3}
F.~Knox and B.~Mohar.
\newblock Maximum number of colourings. {III}. {G}eneral case.
\newblock In preparation.

\bibitem{McKayLink}
B.~McKay.
\newblock Combinatorial data.
\newblock \url{http://users.cecs.anu.edu.au/~bdm/data/ramsey.html}.
\newblock Accessed: 2017-09-22.

\bibitem{RoyleLink}
G.~Royle.
\newblock Colourings of small graphs.
\newblock \url{http://staffhome.ecm.uwa.edu.au/~00013890/remote/graphs/#cols}.
\newblock Accessed: 2017-09-22.

\bibitem{To71}
I.~Tomescu.
\newblock Le nombre des graphes connexes $k$-chromatiques minimaux aux sommets
  \'{e}tiquet\'{e}s.
\newblock {\em C. R. Acad. Sci. Paris}, 273:1124--1126, 1971.

\bibitem{To_DM72}
I.~Tomescu.
\newblock The maximum number of 3-colorings of a connected graph.
\newblock {\em Discrete Mathematics}, 4(1):351--356, 1972.

\bibitem{To_JGT90}
I.~Tomescu.
\newblock Maximal chromatic polynomials of connected planar graphs.
\newblock {\em Journal of Graph Theory}, 14(1):101--110, 1990.

\end{thebibliography}

\newpage
\appendix

\section{Data tables} \label{sec:data tables}

\begin{table}[h]
\centering
\begin{tabular}{|c||l|}
  \hline \\[-0.9em]
  $G[S]$ & $R(y) = \sum_{T \subseteq S'} E'_T y^{4-\Delta'(T)}$ \\[1mm]
  \hline \\[-0.9em]
$5K_1$ & $y^8 - 5y^7 + 10y^6 - 10y^5 + 10y^4 - 11y^3 + 10y^2 - 5y + 1$ \\
$K_2 \cup 3K_1$ & $y^8 - 6y^7 + 16y^6 - 20y^5 + 9y^4 + 5y^3 - 4y^2 - y + 1$ \\
$P_3 \cup 2K_1$ & $y^8 - 7y^7 + 24y^6 - 43y^5 + 38y^4 - 10y^3 - 3y^2 + 1$ \\
$K_3 \cup 2K_1$ & $y^8 - 8y^7 + 33y^6 - 78y^5 + 102y^4 - 71y^3 + 29y^2 - 13y + 4$ \\
$K_{1, 3} \cup K_1$ & $y^8 - 7y^7 + 24y^6 - 41y^5 + 26y^4 + 17y^3 - 29y^2 + 9y + 1$ \\
$K_{1, 4}$ & $y^8 - 7y^7 + 24y^6 - 40y^5 + 22y^4 + 26y^3 - 39y^2 + 14y$ \\
$2K_2 \cup K_1$ & $y^8 - 7y^7 + 23y^6 - 37y^5 + 23y^4 + 18y^3 - 34y^2 + 13y + 1$ \\
$P_4 \cup K_1$ & $y^8 - 8y^7 + 33y^6 - 75y^5 + 93y^4 - 49y^3 - y^2 + 5y + 2$ \\
$K_3 + \text{ leaf } \cup K_1$ & $y^8 - 8y^7 + 33y^6 - 80y^5 + 111y^4 - 77y^3 + 21y^2 - 6y + 5$ \\
$C_4 \cup K_1$ & $y^8 - 9y^7 + 44y^6 - 121y^5 + 190y^4 - 155y^3 + 57y^2 - 11y + 5$ \\
$\overline{K_{1, 4} + e}$ & $y^8 - 8y^7 + 33y^6 - 80y^5 + 113y^4 - 86y^3 + 36y^2 - 18y + 10$ \\
$K_4 \cup K_1$ & $y^8 - 7y^7 + 23y^6 - 45y^5 + 50y^4 - 27y^3 + 17y^2 - 28y + 21$ \\
$P_3 \cup K_2$ & $y^8 - 8y^7 + 32y^6 - 69y^5 + 80y^4 - 28y^3 - 30y^2 + 23y$ \\
$\text{fork}$ & $y^8 - 8y^7 + 33y^6 - 75y^5 + 91y^4 - 32y^3 - 37y^2 + 28y$ \\
$K_{1, 4} + e$ & $y^8 - 8y^7 + 33y^6 - 81y^5 + 118y^4 - 87y^3 + 16y^2 + 8y$ \\
$P_5$ & $y^8 - 9y^7 + 43y^6 - 117y^5 + 186y^4 - 148y^3 + 29y^2 + 16y$ \\
$\text{bull}$ & $y^8 - 8y^7 + 33y^6 - 80y^5 + 110y^4 - 64y^3 - 13y^2 + 21y$ \\
$C_4 + \text{ leaf}$ & $y^8 - 9y^7 + 44y^6 - 124y^5 + 207y^4 - 179y^3 + 50y^2 + 11y$ \\
$\overline{K_3 + \text{ leaf } \cup K_1}$ & $y^8 - 8y^7 + 33y^6 - 80y^5 + 114y^4 - 84y^3 + 18y^2 + 6y$ \\
$K_{2, 3}$ & $y^8 - 9y^7 + 45y^6 - 129y^5 + 219y^4 - 199y^3 + 73y^2$ \\
$\overline{K_3 \cup 2K_1}$ & $y^8 - 8y^7 + 34y^6 - 87y^5 + 137y^4 - 121y^3 + 43y^2$ \\
$C_5$ & $y^8 - 10y^7 + 55y^6 - 180y^5 + 374y^4 - 512y^3 + 415y^2$ \\
$\overline{P_5}$ & $y^8 - 9y^7 + 44y^6 - 130y^5 + 233y^4 - 227y^3 + 85y^2$ \\
$\overline{P_4 \cup K_1}$ & $y^8 - 8y^7 + 33y^6 - 83y^5 + 127y^4 - 106y^3 + 36y^2$ \\
$\overline{P_3 \cup K_2}$ & $y^8 - 8y^7 + 34y^6 - 85y^5 + 123y^4 - 90y^3 + 21y^2$ \\
$\overline{P_3 \cup 2K_1}$ & $y^8 - 7y^7 + 24y^6 - 52y^5 + 66y^4 - 36y^3 - y^2$ \\
$\overline{2K_2 \cup K_1}$ & $y^8 - 7y^7 + 24y^6 - 45y^5 + 44y^4 - 17y^3$ \\
$K_5 - e$ & $y^8 - 6y^7 + 17y^6 - 36y^5 + 57y^4 - 42y^3$ \\
$K_5$ & $y^8 - 5y^7 + 15y^6 - 55y^5 + 109y^4 - 75y^3$ \\
$K_3 \cup K_2$ & $y^8 - 9y^7 + 42y^6 - 114y^5 + 183y^4 - 157y^3 + 41y^2 + 12y$ \\
$\overline{C_4 + \text{ leaf}}$ & $y^8 - 9y^7 + 43y^6 - 123y^5 + 213y^4 - 198y^3 + 60y^2 + 12y$ \\
$\overline{\text{fork}}$ & $y^8 - 8y^7 + 33y^6 - 81y^5 + 117y^4 - 81y^3 + 3y^2 + 15y$ \\
$\overline{K_{1, 3} \cup K_1}$ & $y^8 - 7y^7 + 23y^6 - 45y^5 + 48y^4 - 12y^3 - 23y^2 + 14y$ \\
$\overline{C_4 \cup K_1}$ & $y^8 - 9y^7 + 43y^6 - 125y^5 + 222y^4 - 217y^3 + 85y^2$ \\
  \hline
\end{tabular}
\vskip 2mm
\caption{Polynomials $R(y)$ for each of the 34 graphs on $5$ vertices.}
\label{tab:1}
\end{table}

\begin{table}[h]
\centering
\begin{tabular}{|c||l|}
  \hline \\[-0.9em]
  $G[S]$ & $W(z) = (z+4)^7 (z+1) - R(z + 4)$ \\[1mm] \hline \\[-0.9em]
$5K_1$ & $2z^7 + 46z^6 + 442z^5 + 2270z^4 + 6571z^3 + 10170z^2 + 6597z + 51$ \\
$K_2 \cup 3K_1$ & $3z^7 + 68z^6 + 644z^5 + 3271z^4 + 9451z^3 + 14952z^2 + 10801z + 1539$ \\
$P_3 \cup 2K_1$ & $4z^7 + 88z^6 + 811z^5 + 4022z^4 + 11402z^3 + 17851z^2 + 13048z + 2223$ \\
$K_3 \cup 2K_1$ & $5z^7 + 107z^6 + 966z^5 + 4738z^4 + 13479z^3 + 21751z^2 + 17525z + 4640$ \\
$K_{1, 3} \cup K_1$ & $4z^7 + 88z^6 + 809z^5 + 3994z^4 + 11247z^3 + 17425z^2 + 12463z + 1899$ \\
$K_{1, 4}$ & $4z^7 + 88z^6 + 808z^5 + 3978z^4 + 11142z^3 + 17071z^2 + 11850z + 1464$ \\
$2K_2 \cup K_1$ & $4z^7 + 89z^6 + 829z^5 + 4157z^4 + 11934z^3 + 18986z^2 + 14243z + 2667$ \\
$P_4 \cup K_1$ & $5z^7 + 107z^6 + 963z^5 + 4687z^4 + 13121z^3 + 20461z^2 + 15155z + 2874$ \\
$K_3 + \text{ leaf } \cup K_1$ & $5z^7 + 107z^6 + 968z^5 + 4769z^4 + 13661z^3 + 22247z^2 + 18126z + 4867$ \\
$C_4 \cup K_1$ & $6z^7 + 124z^6 + 1081z^5 + 5110z^4 + 13915z^3 + 21067z^2 + 14931z + 2391$ \\
$\overline{K_{1, 4} + e}$ & $5z^7 + 107z^6 + 968z^5 + 4767z^4 + 13638z^3 + 22148z^2 + 17938z + 4734$ \\
$K_4 \cup K_1$ & $4z^7 + 89z^6 + 837z^5 + 4290z^4 + 12827z^3 + 22003z^2 + 19364z + 6155$ \\
$P_3 \cup K_2$ & $5z^7 + 108z^6 + 981z^5 + 4820z^4 + 13628z^3 + 21486z^2 + 16153z + 3204$ \\
$\text{fork}$ & $5z^7 + 107z^6 + 963z^5 + 4689z^4 + 13136z^3 + 20485z^2 + 15116z + 2784$ \\
$K_{1, 4} + e$ & $5z^7 + 107z^6 + 969z^5 + 4782z^4 + 13719z^3 + 22340z^2 + 18120z + 4768$ \\
$P_5$ & $6z^7 + 125z^6 + 1101z^5 + 5274z^4 + 14612z^3 + 22675z^2 + 16840z + 3312$ \\
$\text{bull}$ & $5z^7 + 107z^6 + 968z^5 + 4770z^4 + 13664z^3 + 22221z^2 + 18003z + 4732$ \\
$C_4 + \text{ leaf}$ & $6z^7 + 124z^6 + 1084z^5 + 5153z^4 + 14147z^3 + 21650z^2 + 15605z + 2676$ \\
$\overline{K_3 + \text{ leaf } \cup K_1}$ & $5z^7 + 107z^6 + 968z^5 + 4766z^4 + 13620z^3 + 22046z^2 + 17706z + 4552$ \\
$K_{2, 3}$ & $6z^7 + 123z^6 + 1065z^5 + 5001z^4 + 13495z^3 + 20075z^2 + 13576z + 1584$ \\
$\overline{K_3 \cup 2K_1}$ & $5z^7 + 106z^6 + 951z^5 + 4643z^4 + 13129z^3 + 20897z^2 + 16216z + 3728$ \\
$C_5$ & $7z^7 + 141z^6 + 1212z^5 + 5706z^4 + 15648z^3 + 24353z^2 + 18696z + 4112$ \\
$\overline{P_5}$ & $6z^7 + 124z^6 + 1090z^5 + 5247z^4 + 14739z^3 + 23535z^2 + 18664z + 4720$ \\
$\overline{P_4 \cup K_1}$ & $5z^7 + 107z^6 + 971z^5 + 4813z^4 + 13914z^3 + 22964z^2 + 19136z + 5440$ \\
$\overline{P_3 \cup K_2}$ & $5z^7 + 106z^6 + 949z^5 + 4617z^4 + 13002z^3 + 20611z^2 + 15928z + 3632$ \\
$\overline{P_3 \cup 2K_1}$ & $4z^7 + 88z^6 + 820z^5 + 4174z^4 + 12420z^3 + 21233z^2 + 18632z + 5904$ \\
$\overline{2K_2 \cup K_1}$ & $4z^7 + 88z^6 + 813z^5 + 4056z^4 + 11633z^3 + 18636z^2 + 14384z + 3136$ \\
$K_5 - e$ & $3z^7 + 67z^6 + 636z^5 + 3303z^4 + 10010z^3 + 17304z^2 + 15072z + 4480$ \\
$K_5$ & $2z^7 + 41z^6 + 367z^5 + 1871z^4 + 5851z^3 + 11044z^2 + 11280z + 4544$ \\
$K_3 \cup K_2$ & $6z^7 + 126z^6 + 1122z^5 + 5457z^4 + 15469z^3 + 24979z^2 + 20252z + 5504$ \\
$\overline{C_4 + \text{ leaf}}$ & $6z^7 + 125z^6 + 1107z^5 + 5367z^4 + 15190z^3 + 24492z^2 + 19764z + 5264$ \\
$\overline{\text{fork}}$ & $5z^7 + 107z^6 + 969z^5 + 4783z^4 + 13729z^3 + 22377z^2 + 18185z + 4820$ \\
$\overline{K_{1, 3} \cup K_1}$ & $4z^7 + 89z^6 + 837z^5 + 4292z^4 + 12844z^3 + 22055z^2 + 19434z + 6200$ \\
$\overline{C_4 \cup K_1}$ & $6z^7 + 125z^6 + 1109z^5 + 5398z^4 + 15385z^3 + 25111z^2 + 20744z + 5872$ \\
  \hline
\end{tabular}
\vskip 2mm
\caption{All coefficients of $W(z)$ are non-negative, thus showing that $R(y) \ll y^7 (y-3)$ for each graph in the table.}
\label{tab:2}
\end{table}

\end{document}